\documentclass{amsart}                          

\usepackage{amsmath, amscd, amsfonts, amssymb}
\usepackage{mathrsfs}
\usepackage{bbm}
\usepackage{stmaryrd}
\usepackage{graphicx}
\usepackage[all]{xy}
\usepackage{color}

\renewcommand{\arraystretch}{1.3}

\newtheorem{prop}{Proposition}[section]
\newtheorem{thm}[prop]{Theorem}
\newtheorem{lem}[prop]{Lemma}
\theoremstyle{definition}
\newtheorem{defn}[prop]{Definition}

\theoremstyle{remark}
\newtheorem{rem}[prop]{Remark}

\newtheorem{exmp}[prop]{Example}

\newcommand{\z}{\mathbb Z}
\renewcommand{\r}{\mathbb R}
\newcommand{\n}{\mathbb N}
\renewcommand{\c}{\mathbb C}

\renewcommand{\o}{\mathfrak{o}}
\newcommand{\e}{\mathfrak{e}}
\newcommand{\X}{\mathfrak{X}}
\newcommand{\red}{{\rm red}}
\newcommand{\defo}{{\rm def}}
\newcommand{\q}{{\rm q}}

\begin{document}

\title{SUSY Structures on Deformed Supermanifolds}
\author{Frank Klinker}
\thanks{Frank Klinker, Department of Mathematics, University of Dortmund, D--44221 Dortmund}
\date{December 31, 2007}
\address{Department of Mathematics, University of Dortmund, D--44221 Dortmund} 
\address{frank.klinker@mathematik.uni-dortmund.de}
\begin{abstract}
We construct  a geometric structure  on  deformed supermanifolds  as a certain subalgebra of the vector fields.   
In the classical limit we obtain a decoupling of the infinitesimal odd and even transformations, 
whereas in the  semiclassical limit the result is a representation of the supersymmetry algebra. 
In the case that the structure is mass preserving we describe all high energy corrections to this algebra.
\end{abstract}
\maketitle

\section{Introduction and Preliminaries}\label{intro}

In this text we present the definition of a geometric structure on deformations of supermanifolds which we call SUSY structure. More precisely, we construct a subalgebra $\mathfrak{S}$ of the vector fields on the deformed supermanifold, such that in the semiclassical limit we obtain a representation of the supersymmetry algebra. We develop in detail the tools we need to formulate the construction of $\mathfrak{S}$ before we discuss examples. In particular, as an important example of such structure -- a structure which preserves  mass dimension -- we will describe all high energy corrections in detail. The physical notations which we used above are naturally motivated after introducing the mass dimension in section \ref{deformation}.  In particular the different limits  are closely related to certain ranges of values of the deformation parameter which may be interpreted as energy. Their meaning  will be specified in definition \ref{classical}.

Consider the graded manifold $M=(M_\red,\mathfrak{A})$ where $M_\red$ is a (pseudo) Riemannian spin manifold. By $\nabla$ we denote the Levi-Civita connection on $M_\red$. Furthermore $\mathfrak{A}:=\Gamma\Lambda S$ is the space of sections in the exterior bundle of the spinor bundle $S$ over $M_\red$. The splitting of $\mathfrak{A}=\mathfrak{A}_0\oplus \mathfrak{A}_1$ into even and odd forms define the even and odd functions on $M$. The derivations of $\mathfrak{A}$ are called vector fields on $M$ and will be denoted by $\X(M)$. 
The vector fields  inherit a natural $\z_2$-grading $\X(M)=\X(M)_0\oplus \X(M)_1$ which is given by  derivations which preserve  or change the degree  of homogeneous elements in $\mathfrak{A}$, respectively.

The graded manifold $M$ is equipped with a bilinear form $g+C$ where $g$ is the metric on $M$ and $C= \langle\cdot,\cdot\rangle$  a charge conjugation on $S$. 
The charge conjugation is a non degenerated, spin invariant, bilinear map which yields an identification $\Gamma S^*\simeq \Gamma S$. 
It provides a morphism $C_*:\Gamma S\otimes \Gamma S\simeq \Gamma  End(S) \to \Gamma \Lambda^*TM_\red$ and we denote the projections on the respective summands by $C_k:\Gamma S\otimes \Gamma S\to \Gamma \Lambda^kTM_\red$. The induced projections  are either symmetric or skew symmetric. This depends on the choice of charge conjugation and we denote the symmetry by $\Delta_k\in\{\pm1\}$. 
The morphism $C_*$ needs a few words of explanation.
If the spinors are taken to be complex, the target space of $C_*$ should of course be taken complexified, too. In this case, and if $S$ is the minimal representation, $C_*$ is an isomorphism if $\text{dim} M$ is even but maps onto one half of the forms if $\text{dim} M$ is odd -- see \cite{LawMich} or \cite{Chevalley} for more details. 
If $S$ is taken to be real the discussion is more delicate. Depending on the signature of the metric the spin representation is real, complex or quaternionic and a construction of real spinors, in particular in the complex and quaternionic case, needs a doubling of the complex spinors. In these cases $C_*$ is not 1:1 but 1:$2^\alpha$ with $\alpha=0,1,2$ or $3$. For an elaborated discussion  on the restriction to real spinors as well as the symmetry of the projections  we cordially refer to \cite{AlekCortDevProey} and \cite{Klinker4}. 

In the case $\Delta_1=1$ we will often write $C_1(\eta\otimes\xi)=\frac{1}{2}\{\eta,\xi\}$ and we will call $M$ a {\em special graded manifold}. For example in  Lorentzian signature $\Delta_1=1$ is always possible for the minimal real spinors, see \cite{AlekCor2} and \cite{Klinker4}. 

A map  $\jmath_C:\Gamma S\to\X(M)_1$ is now canonically assigned to $C_0=C$. More precisely, $\jmath_C(\varphi)$ acts via interior multiplication with respect to the charge conjugation and, therefore, is a derivation of degree $-1$.

Consider a connection $D$ on the spinor bundle $S$. For every vector field $X\in\X(M_\red)$ the action of $D_X$ on $\mathfrak{A}$ is that of a derivation of degree zero. We  write $\jmath_D:\X(M_\red) \to\X(M)_0$ for the corresponding map.

The inclusion $\jmath_D\oplus\jmath_C:\X(M_\red) \oplus \Gamma S\hookrightarrow \X(M)$ gives rise to a splitting\footnote{
We will omit the index at the inclusion maps if we fix the charge conjugation and the connection.}
\begin{equation}\label{splitting}
\X(M) = \mathfrak{A} \otimes \jmath \X(M_\red)  \oplus \mathfrak{A} \otimes \jmath \Gamma S\,,
\end{equation}
see \cite{Kostant} or \cite{ManSardan}.
We call a vector field $X$ of degree $(k,\epsilon)\in\z\times\z_2$ 
if $X\in\Gamma\Lambda^kS\otimes \Gamma S$ for $\epsilon=1$ 
and $ X\in\Gamma\Lambda^k S\otimes \X(M_\red) $ for $\epsilon=0$. 
The $\z$-degree of a vector field of degree $(k,\epsilon)$ is defined as $k-\epsilon\in\z$, 
and  $(k-\epsilon)\,\text{mod}\,2$ coincides with the induced $\z_2$-grading. 
The  even and odd parts of the vector fields are given by
\begin{equation}
\X(M)_{0/1} = 		\mathfrak{A}_{0/1} \otimes \X(M_\red) 
		\oplus 	\mathfrak{A}_{1/0}\otimes \Gamma S\,,
\end{equation}
whereas the vector-  and  spinor-like fields are defined as $
\X(M)^{\text{v}}:=\mathfrak{A}\otimes\X(M_\red)$ and $\X(M)^{\text{s}}:=\mathfrak{A}\otimes \Gamma S$.

\section{Filtrations and Deformations}\label{deformation}

{\bf Filtrations.\ }\ 
A {\em filtration} of a set $X$ is an increasing sequence $\mathcal{X}=(X_n)_{n\in\n}$ of subsets of $X$, i.e.\  $X_i\subset X_{i+1}$, such that $\bigcup_kX_k=X$. A filtration is called {\em finite} if there exits a $n_0$ such that $X_k=X$ for all $k\geq n_0$.
A filtration of a linear space $X$ is an increasing sequence of linear subspaces which, in particular, is finite if $X$ is finite dimensional. 
If $X$ has an algebra structure the filtration is called $\ell$-{\em compatible} if $\ell$ is the minimal number such that 
$X_{i}X_{i'}\subset X_{i+i'+\ell}$ for all $i,i'\in\n$. 
A subsequence $\mathcal{X}'=(X'_n)_{n\in\n}$  of  $\mathcal{X}$ is called {\em subfiltration}. 
If the  filtrations $\mathcal{X}=(X_n)_{n\in\n}$ and $\mathcal{X}'=(X'_n)_{n\in\n}$ of $X$ obey $X_i\subset X'_{i}\subset X_{i+1}$ for all $i\in\n$ we say $\mathcal{X'}$ {\em supplements} $\mathcal{X}$. 
If $X$ admits a $\z_2$-grading we call a filtration $\mathcal{X}=\{X_n\}_{n\in\n}$ {\em compatible with the grading} if $X_k^{k\,\text{mod}\,2}=X_{k+1}^{k\,\text{mod}\,2}$ for all $k$.

We will apply these notions to the graded manifold $M$ and discuss filtrations of the algebra of functions as well as of the super Lie algebra of vector fields.

Consider the graded manifold $M=(M_\red,\mathcal{A})$. The algebra of superfunctions  admits a natural $0$-compatible filtration $\mathfrak{A}=\bigoplus_k \mathfrak{A}_k$ with 
\begin{equation}
\mathfrak{A}_k:=\bigoplus\nolimits_{\ell\leq k}\Gamma\Lambda^\ell S\,.\label{natural}
\end{equation}
The Lie algebra of vector fields $\X(M)$ on the graded manifold admits different filtrations.
We consider two filtrations given by $\mathcal{W}=\{W_i\}_{i\geq0}$ and $\widehat{\mathcal{W}}=\{\widehat{W}_i\}_{i\geq-1}$ with\footnote{We write $\X:=\X(M_\red)$ in this subsection.}  
\begin{equation*}
{\setlength{\arraycolsep}{0.5ex}\begin{array}{rrcr}
W_{4k}		:=& \mathfrak{A}_{2k}\otimes\jmath(\X ) 
			&\oplus&\mathfrak{A}_{2k-1}\otimes\jmath(\Gamma S)\,, \\
W_{4k+1}	:=& \mathfrak{A}_{2k}\otimes\jmath(\X ) 
			&\oplus&\mathfrak{A}_{2k}\otimes\jmath(\Gamma S)\,, \\
W_{4k+2}	:=& \mathfrak{A}_{2k}\otimes\jmath(\X ) 
			&\!\!\!\oplus\!\!\!&\mathfrak{A}_{2k+1}\otimes\jmath(\Gamma S)\,, \\
W_{4k+3}	:=& \mathfrak{A}_{2k+1}\otimes\jmath(\X ) 
			&\oplus&\mathfrak{A}_{2k+1}\otimes\jmath(\Gamma S)\,, \\
\end{array}}
\end{equation*}
and
\begin{equation*}
{\setlength{\arraycolsep}{0.5ex}\begin{array}{rrcr}
\widehat{W}_{4k-1}:=& \mathfrak{A}_{2k-1}\otimes\jmath(\X) 
			&\oplus& \mathfrak{A}_{2k}\otimes\jmath(\Gamma S )\,, \\
\widehat{W}_{4k}	:=&  \mathfrak{A}_{2k}\otimes\jmath(\X)
			&\oplus&\mathfrak{A}_{2k}\otimes\jmath(\Gamma S )\,,  \\
\widehat{W}_{4k+1}:=& \mathfrak{A}_{2k+1}\otimes\jmath(\X)
			&\oplus&  \mathfrak{A}_{2k}\otimes\jmath(\Gamma S )\,, \\
\widehat{W}_{4k+2}:=&\mathfrak{A}_{2k+1}\otimes\jmath(\X)  
			&\oplus& \mathfrak{A}_{2k+1}\otimes\jmath(\Gamma S )\,.
\end{array}}
\end{equation*}
The filtrations $\mathcal{W}$ and $\widehat{\mathcal{W}}$ have the following properties:
\begin{itemize}
\item $\mathcal{W}$ and $\widehat{\mathcal{W}}$ are finite with $\X(M)=W_{2\ell+1}=\widehat{W}_{2\ell}$ for $\ell\geq\dim S$.
\item Both filtrations are fine in the sense that $W_{m+1}/W_m$ is homogeneous in $\X(M)$ with respect to the $(\z\times\z_2)$-grading.
\item Both filtrations are compatible with the $\z_2$-structure, i.e.\ $W_{2k}^0=W_{2k+1}^0$ and $W_{2k-1}^1=W_{2k}^1$ as well as
$\widehat{W}_{2k}^0=\widehat{W}_{2k+1}^0$ and $\widehat{W}_{2k-1}^1=\widehat{W}_{2k}^1$
\item
With respect to the graded Lie bracket --  see (\ref{basiccomm}) for some basic formulas -- 
$\mathcal{W}$ is $2$-compatible and  the two subfiltrations $\{W_{2n+1}\}$ and $\{W_{2n}\}$ are $1$-compatible.  $\widehat{\mathcal{W}}$ is $5$-compatible, although its subfiltration $\{\widehat{W}_{2n}\}$  $(=\{W_{2n+1}\})$  is $1$-compatible.
\end{itemize}
The last point makes the filtration $\mathcal{W}$ the more natural one. Nevertheless, we will need a combination of both.

By defining
$G_i:=W_{2i}^0=\widehat{W}_{2i}^0$
and 
$U_i:=W_{2i+1}^1=\widehat{W}_{2i-1}^1$ 
we get an even filtration $\{G_i\}$  of $\X(M)_0$ and an odd filtration $\{U_i\}$ of $\X(M)_1$ which are explicitly given by 
\begin{equation*}
{\setlength{\arraycolsep}{0.5ex}\begin{array}{rrcr}
G_{2k}		=&\mathfrak{A}_{2k}^0 \otimes\jmath(\X) 
			&\oplus&\mathfrak{A}_{2k-1}^1\otimes \jmath(\Gamma S)\,, \\
G_{2k+1}	=& \mathfrak{A}^0_{2k}\otimes\jmath(\X) 
			&\oplus& \mathfrak{A}_{2k+1}^1\otimes \jmath(\Gamma S)\,, \\
U_{2k}		=& \mathfrak{A}_{2k-1}^1\otimes\jmath(\X) 
			&\oplus&\mathfrak{A}_{2k}^0 \otimes\jmath(\Gamma S)\,,\\
U_{2k+1}	=&\mathfrak{A}_{2k+1}^1\otimes\jmath(\X)
			&\oplus& \mathfrak{A}_{2k}^0\otimes\jmath(\Gamma S) \,.
\end{array}}
\end{equation*}
The filtrations are connected via $G_k\oplus U_k=W_{2k+1}=\widehat{W}_{2k}$ and $G_{k}\oplus U_{k-1}=W_{2k}$ as well as $ G_{k}\oplus U_{k+1}=\widehat{W}_{2k+1}$. In particular,
\begin{align*}
G_{2k+1}\oplus U_{2k} 	&= W_{4k+2} 
			= \mathfrak{A}_{2k}\otimes\jmath(\X)\oplus \mathfrak{A}_{2k+1}\otimes\jmath(\Gamma S)\,,\\
G_{2k+1}\oplus U_{2k+2}	&=\widehat{W}_{4k+3}
			= \mathfrak{A}_{2k+1}\otimes\jmath(\X)\oplus \mathfrak{A}_{2k+2}\otimes\jmath(\Gamma S)\,.
\end{align*}
\begin{defn}\label{defdef}
We define the filtration $\mathcal{Z}=\{Z_k\}_{k\geq0}$ by
\begin{equation}
Z_{2k}:=W_{4k+2}, \qquad  Z_{2k+1}:=\widehat{W}_{4k+3}\,.
\end{equation}
I.e.\ $Z_k=\mathfrak{A}_k\otimes \jmath(\X)\oplus\mathfrak{A}_{k+1}\otimes\jmath(\Gamma S)$.
\end{defn} 
The filtration $\mathcal{Z}$ has the following properties:
\begin{itemize}
\item
$\mathcal{Z}$ is neither a subfiltration of $\mathcal{W}$ nor $\widehat{\mathcal{W}}$  but a combination of parts of both.
\item 
$\mathcal{Z}$ is compatible with the $\z_2$-grading, i.e.\ $Z_{2k}^0=Z_{2k+1}^0$ and $Z_{2k-1}^1=Z_{2k}^1$.
\item
$\mathcal{Z}$ is not fine with respect to the $(\z\times\z_2)$-grading but it is  fine with respect to the $\z$-grading, i.e.\ the $(\z\times\z_2)$-homogeneous elements in  $Z_{k}/Z_{k-1}$ are of degree $(k,0)$ or $(k+1,1)$.
\item $\mathcal{Z}$ is a $0$-compatible filtration.
\end{itemize}


{\bf Deformations.\ }\  
\begin{defn}
 The {\em deformation} of a set $X$ with respect to the filtration $\mathcal{X}$ is defined as the subset 
\begin{equation}
\defo_\q^\mathcal{X}X:=\big\{{\textstyle \sum} a_k\q^k\,|\,a_k\in X_k \big\}
\end{equation}
of $X\llbracket\q\rrbracket$ of formal power series in the (even)  parameter $\q$ with coefficients in $X$. 
\end{defn}
\begin{rem}\begin{itemize}
\item Deformations as given above have been introduced by M.\ Gerstenhaber in a series of papers of which we would like to mention  the last one \cite{GerstenhaberV}. The author uses  filtrations which obey $F_j\subset F_{j-1}$ with respect to a $\z$-grading. Up to sign this coincides with our $\z$-grading. 
\item
If $X$ is a algebra and $\mathcal{X}=\{X_n\}_{n\in\n}$ is $0$-compatible then $\defo_\q^\mathcal{X}X$ is a subalgebra of $X\llbracket\q\rrbracket$.
\item
If $\mathcal{X}$  is compatible with a given $\z_2$-grading on $X$  then $\defo_\q^\mathcal{X}X$ is canonically $\z_2$-graded.
\end{itemize}
\end{rem}

\begin{defn} 
Let $M=(M_\red,\mathfrak{A})$ be a graded manifold. 
A {\em deformation} of $M$  is defined by the deformation of $\mathfrak{A}$ with respect to the natural  filtration (\ref{natural}). 
We write  $M_\q:=(M_\red, \defo_\q\mathfrak{A})$  with
\begin{equation}
\defo_\q \mathfrak{A}:=\defo_\q^{\text{nat}}\mathfrak{A} \subset \mathfrak{A}\llbracket\q\rrbracket\ .
\end{equation}
The vector fields on the deformed graded manifold $M_\q$ are  given by the derivations of the algebra $\defo_\q\mathfrak{A}$, and will be denoted by $\X(M_\q)= \mathfrak{der}\big( \defo_\q \mathfrak{A}\big)$.
\end{defn}

\begin{prop}
Let $M_\q$ be the deformation of the graded manifold $M=(M_\red,\mathfrak{A})$ and $\mathfrak{A}=\Gamma\Lambda S$. 
The vector fields of the deformed graded manifold $M_\q$ are connected to the filtration $\mathcal{Z}$ of $\X(M)$ cf.\ definition \ref{defdef} by
\begin{equation}
\X(M_\q)  \cap \X(M)\llbracket\q\rrbracket = \defo_\q^\mathcal{Z}\X(M) 
\end{equation}
\end{prop}

\begin{proof}
The deformation of $\X(M)$ with respect to $\mathcal{Z}$ is given by 
\begin{align*}
&\defo_\q^\mathcal{Z}\X(M) \\
	&\ = \big\{ {\textstyle\sum} a_k\q^k\,|\,a_k\in Z_k\}\\
	&\ = \Big\{{\textstyle \sum} \big(a_k\otimes \jmath(X)+b_{k+1}\otimes\jmath(\theta)\big)\q^k\,|\,
			a_\ell,b_\ell\in\mathfrak{A}_\ell,X\in\X(M_\red),\theta\in \Gamma S\Big\}\\
	&\ = \Big\{ {\textstyle \sum}  a_k\q^k\otimes \jmath(X)+{\textstyle \sum} b_{k+1}\q^k\otimes\jmath(\theta)\,|\,
			a_\ell,b_\ell\in\mathfrak{A}_\ell , X\in\X(M_\red),\theta\in \Gamma S\Big\}\\
	&\ = \Big\{ a \otimes \jmath(X)+b\otimes\jmath(\theta)\,|\,
			a\in\defo_\q\mathfrak{A},b\in\mathfrak{A}_1\cdot\defo_\q\mathfrak{A},
			 X\in\X(M_\red),\theta\in \Gamma S\Big\}\\
	&\ = \defo_\q\mathfrak{A}\otimes \jmath(\X(M_\red)) 
		\ \oplus \ \mathfrak{A}_1\cdot\defo_\q\mathfrak{A}\otimes\jmath(\Gamma S)\,.
\end{align*}
The derivations of $\defo_\q \mathfrak{A}$ form a subset of $\mathfrak{der}\big(\mathfrak{A}\llbracket\q\rrbracket\big)$. A careful examination of the possible homogeneous summands leads to
\begin{equation*}
{\setlength{\arraycolsep}{0.5ex}\begin{array}{rccc}
\mathfrak{der}\big( \defo_\q \mathfrak{A}\big)
	=&\underbrace{\defo_\q\mathfrak{A}\otimes \jmath(\X(M_\red)) 
		\ \oplus \ \mathfrak{A}_1\cdot\defo_\q\mathfrak{A}\otimes\jmath(\Gamma S)}
		&\oplus& \q\cdot\defo_\q\mathfrak{A}\cdot \frac{\partial}{\partial\q}\\
&\cap &&\cap\\
\mathfrak{der}\big(\mathfrak{A}\llbracket\q\rrbracket\big)
	=& \X(M) \llbracket\q\rrbracket 
		&\oplus& \mathfrak{A}\llbracket\q\rrbracket\, \frac{\partial}{\partial\q}
\end{array}}
\end{equation*}
Comparing the two spaces finishes the proof.
\end{proof}


{\bf Mass dimension.\ }\  It is natural to express the units of all physical values in powers  of the unit of mass and we call this power the {\em mass dimension}. 
For example, we have $[\text{energy}]=[\text{momentum}]=[\text{acceleration}]=+1$ or $[\text{time}]=[\text{length}]=-1$ or $[\text{velocity]}=0$. In this context a vector field as a field of directional derivatives or velocities  has mass dimension $0$. 
It is natural to assign to a spinor  mass dimension $-\frac{1}{2}$ so that a a supersymmetry generator has mass dimension $+\frac{1}{2}$.
A multi spinor $\theta\in\Lambda^k\Gamma S$, thus, has  mass dimension $-\frac{k}{2}$. 

The even formal parameter $\q$ is defined to be of mass dimension $[\q]:=1$. 
With this notation we get the following possible mass dimensions for vector fields of $M_\q$ contained in $\defo_\q^\mathcal{Z}\X(M)$:
\begin{table}[htb]\caption{Mass dimensions}\label{tablemass}
$\begin{array}{rcl}
\displaystyle \theta\in\q^{2\ell}Z^0_{2\ell}	& \Rightarrow& \displaystyle [\theta]\in\big\{\ell,\ell+1,\ldots,2\ell\big\}\\
\displaystyle \theta\in\q^{2\ell}Z^1_{2\ell}	& \Rightarrow& \displaystyle [\theta]\in\big\{\ell+\tfrac{1}{2},\ell+\tfrac{3}{2},\ldots,2\ell+\tfrac{1}{2}\big\}\\
\displaystyle \theta\in\q^{2\ell+1}Z^0_{2\ell+1}	& \Rightarrow&\displaystyle [\theta]\in\big\{\ell+1,\ell+2,\ldots,2\ell+1\big\}\\
\displaystyle \theta\in\q^{2\ell+1}Z^1_{2\ell+1}	& \Rightarrow&\displaystyle  [\theta]\in\big\{\ell+\tfrac{1}{2},\ell+\tfrac{3}{2},\ldots, 2\ell+\tfrac{3}{2}\big\}
\end{array}$
\end{table}
\begin{rem} 
The mass dimension of vector fields in $\X(M)$ coincides with the $\z$-degree up to a factor $-2$. So the  mass-dimension coincides with grading cf.\ Gerstenhaber, at least up to a factor of $2$.
Therefore, the mass dimension on $\defo_\q^\mathcal{Z}\X(M)$  may be identified with an extension of the $\z$-grading on $\X(M)$ where $\q$ has degree $-2$. This extension is compatible with the $\z_2$-grading because $\q$ is even.
\end{rem}
Because $\q$ has mass dimension one it may be identified with the physical value energy. This motivates the following general definition.
\begin{defn}\label{classical}
For a subset $U\subset \defo_\q^\mathcal{X}X$ of a deformation of $X$ with respect to a filtration $\mathcal{X}$ we call
$U\,\text{mod}\,\q$  the {\em classical limit} and 
$U\,\text{mod}\,\q^2$ the {\em semiclassical limit}.
The terms which are proportional to $\q^k$ with $k\geq2$ are called {\em high energy contributions}.  If we want to emphasize that the high energy contributions turn $U$  into an algebra we  call them {\em high energy corrections}. 
\end{defn}


\section{Torsion and Admissibility}

In this section we will give the construction of torsion of arbitrary spinor connections. The result of the construction is an extension of the notion of torsion for connections on the tangent bundle of a pseudo Riemannian manifold. It is a proper extension in such a way that for metric connections on the spinor bundle, i.e.\ those which are lifted from a connection on $TM$ via the Clifford map, both notions coincide. The torsion of an arbitrary connection $D$ depends on the choice of a spin invariant bilinear form $C$ on the spinor bundle (see definition \ref{ctorsion}). This is the price  for the  fact that the  class of connections we will consider does not preserve the charge conjugation in general, see definition \ref{definitionadmissible}, remark \ref{remarkadmissible}, and   proposition \ref{propositionadmissible}. 
For the proofs of the statements and an extended discussion of the objects which we recall below we cordially refer the reader to \cite{Klinker5} where the torsion of spinor connections has been introduced.

As in section \ref{intro} we denote the  charge conjugation on the spinor bundle $S$  by $C$.  Moreover, we take the  Clifford map as $\gamma\in TM\otimes S\to S$ or $\gamma:TM\to End S$ depending on what is more convenient in the respective situation. For the images of a local base $\{e_k\}$ of $TM$ we use the common notation $\gamma(e_k)=\gamma_k$. 

Given a spinor connection on $S$, we denote the  connection which is  induced by $D$ and the Levi-Civita connection on bundles constructed of $S$ and $TM_\red$ by algebraic operations by the same symbol $D$, e.g. $D(\eta\otimes X)=D\eta\otimes X+\eta\otimes \nabla X$ for $\eta\in\Gamma S$ and $X\in\X(M_\red)$ or $(D\alpha)(\eta)=-\alpha(D\eta)$ for $\alpha\in\Gamma S^*$ and $\eta\in\Gamma S$.

We will assign  to the connection $D$ another connection $D^C$ on $S$ by requiring that the connection $D\otimes \mathbbm{1}+\mathbbm{1}\otimes D^C$
on $S\otimes S$ makes the charge conjugation parallel, i.e.\ $C(D_X\eta,\xi)+C(\eta,D^C_X\xi)=0$ for all $\eta,\xi\in\Gamma S$ and $X\in\X(M_\red)$. 
The next remark is obtained immediately.

\begin{rem}\label{curvC}
For $\Omega\in \Gamma End(S)$ we define  $\Omega^C\in \Gamma End(S)$ by $C(\Omega^C\eta,\xi):=C(\eta,\Omega\xi)$ for all $\xi,\eta\in\Gamma S$. 
Then the curvature $R$ of  $D$ and the curvature $R^C$ of  $D^C$ are related by 
\[
(R(X,Y))^C=-R^C(X,Y)\,.
\] 
\end{rem}

We endow the bundle of ${End}(S)$-valued tensors with a connection $\hat D$ which is induced by $D$, $D^C$ and $\nabla$ in the following way.
\begin{defn}
Let $\Phi\in \X(M_\red)^{\otimes k} \otimes \Omega^1(M_\red)^{\otimes\ell} \otimes {End}(S)$. The connection $\hat D$ is defined by 
\begin{equation}
(\hat D_Z\Phi)(X)\xi := D_Z(\Phi(X)\xi)-\Phi(\nabla^0_ZX)\xi-\Phi(X)D^C_Z\xi 
\label{definitionhatD}
\end{equation}
for all spinors $\xi\in S$,  tensor fields $X\in \Omega^1(M_\red)^{\otimes k}\otimes \X(M_\red) ^{\otimes\ell}$,  and vector fields $Z\in\X(M_\red)$.
\end{defn}
We consider the representation $ad^C$ of ${End}(S)$ on itself given by 
\begin{equation}
ad^C_\Omega: {End}(S)\to {End}(S),\qquad
ad^C_\Omega \Phi:= \Omega\Phi+\Phi \Omega^C\,.
\end{equation}

\begin{rem}
\begin{itemize}
\item
$ad^C$ and $\hat D$ are compatible in the following way:
\begin{equation}\label{compDad}
\hat D(ad^C_\Omega\Psi)=ad^C_{D\Omega}\Psi+ad^C_\Omega \hat D\Psi\,.
\end{equation}
\item 
Denote by $\Omega^\pm$ the two respective projections on the $(\pm1)$-eigen\-spaces of $(\cdot)^C$ . Then 
$ad^C_\Omega \Phi= \big[\Omega_-,\Phi\big]+\big\{\Omega_+,\Phi\big\} $. 
\item 
Furthermore we have
$ (ad^C_\Omega \Phi)^C =ad^C_\Omega\Phi^C $ which yields that
$ad^C_\Omega$ preserves the $(\pm 1)$-eigenspaces of the linear map $\Phi\mapsto\Phi^C$ for all $\Omega\in{End}(S)$.
\end{itemize}
\end{rem}

\begin{defn}\label{ctorsion}
Let $D$ be a connection on the spinor bundle $S$ over the manifold $M_\red$. Denote the Levi-Civita connection on $M_\red$ by $\nabla$  and the Clifford map by $\gamma:\X(M_\red) \to\Gamma{End}(S)$.
The {\em torsion} $\mathcal{T}\in\Omega^2(M_\red)\otimes\Gamma {End}(S)$ of  $D$ is defined by two times the skew symmetrization of $\hat D\gamma:\X(M_\red) \otimes \X(M_\red) \to \Gamma{End}(S)$, i.e.
\begin{equation}
\mathcal{T}(X,Y)= (\hat D_X\gamma)(Y)-(\hat D_Y\gamma)(X)\,.
\end{equation}
\end{defn}

\begin{rem}\label{remarktorsion}
\begin{itemize}
\item We write $(\hat D_X\gamma)(Y)=\hat D_X(\gamma(Y)) - \gamma(\nabla_XY)$. If we omit the map $\gamma$ we get
\[
\mathcal{T}(X,Y)=\hat D_XY-\hat D_YX-[X,Y]\,.
\]
\item In terms of the difference $\mathcal{A}=D-\nabla \in\Omega^1(M)\otimes\Gamma{End}(S)$  the torsion may be written as
\[
\mathcal{T}(X,Y)= ad^C_{\mathcal{A}(X)}Y-ad^C_{\mathcal{A}(Y)}X\,.
\]
\item
The torsion has symmetry $\Delta_1$, i.e.\ we have  
$ C(\eta,\mathcal{T}_{\mu\nu}\xi)=\Delta_1C(\xi,\mathcal{T}_{\mu\nu}\eta)$ for all $\eta,\xi$.
\item 
For a metric connection $D$ on $S$ the torsion coincides with  the torsion  defined by the connection $D$ on $M_\red$. 
\end{itemize}
\end{rem}

The torsion obeys some  Bianchi-type identities.

\begin{prop}\label{Bianchi}
Let $D$ be a connection on the spinor bundle $S$ over the \linebreak (pseudo) Riemannian manifold $M_\red$. The torsion  $\mathcal{T}$ and the curvature $R$ of $D$ obey
\begin{gather}
\hat D_{[\kappa}\mathcal{T}_{\mu\nu]}= ad^C( R_{[\kappa\mu})\gamma_{\nu]}\,, \label{BianchiDT}\\
\hat D_{[\kappa}(ad^C_ R\gamma)_{\mu\nu\rho]}=ad^C(R_{[\kappa\mu})\mathcal{T}_{\nu\rho]}\label{BianchiDadR}\,.
\end{gather}
In this context we add the following identity for the curvature $R$ of $D$:
\begin{equation}
D_{[\kappa}R_{\mu\nu]}=0\,.\label{BianchiDR}
\end{equation}
\end{prop}

The difference $\mathcal{A}=D-\nabla$ of the spinor connection $D$ and the Levi-Civita  connection $\nabla$ is an $End(S)$-valued one-form on $M_\red$. Therefore $\mathcal{A}(X)$ can be decomposed into its homogeneous summands via the morphism $C_*$. 
Usually, for a given $\ell$-form $F$ on $M_\red$ the two terms $X\wedge F$ and $X\lfloor F$ both contribute to physically relevant connections.\footnote{For example the supergravity connection in eleven dimension has a four-form contribution via $D_\mu-\nabla_\mu=-\tfrac{1}{36}F_{\mu\nu\rho\sigma}\gamma^{\nu\rho\sigma}+\tfrac{1}{288}F^{\nu\rho\sigma\tau}\gamma_{\mu\nu\rho\sigma\tau}$, see \cite{CremmJulScher}.} 
Due to the fact that either $X\wedge F$ or $X\lfloor F$ preserve the charge conjugation -- but in no case both of them -- the class of connections which is adapted to our purpose will be different. 

\begin{defn}\label{definitionadmissible}
Let $D$ be a connection on $S$ and $\mathcal{K}_1\subset \Gamma S$ a subset. 
We call  the pair $(D,\mathcal{K}_1)$  {\em admissible} if the symmetric part of $\hat D\gamma $ acts trivially on $\mathcal{K}_1$ and if $D$ is fixed, we call $\mathcal{K}_1$ admissible. 
We call $D$ {\em admissible} if $\hat D\gamma$ is skew symmetric
\end{defn}

\begin{rem}\label{remarkadmissible}
\begin{itemize}
\item 
If $(D,\mathcal{K}_1)$ is admissible and $\mathcal{T}$ is the torsion of the connection $D=\nabla+\mathcal{A}$ we have 
$\mathcal{T}_{\mu\nu}\xi = 2\hat D_\mu\gamma_\nu \xi = 2 ad^C_{\mathcal{A}_\mu}\gamma_\nu\xi $ for all $\xi\in\mathcal{K}_1$.
\item 
If the connection $D$ on $S$ is admissible so is $(D,\mathcal{K}_1)$ for all $\mathcal{K}_1\subset\Gamma S$.
\item 
Admissible metric connections are exactly those metric connections with totally skew symmetric torsion.
\end{itemize}
\end{rem}
In particular, the next two propositions makes  admissibility appropriate concept for our main purpose.  
\begin{prop}
The connection $D=\nabla+\mathcal{A}$ on $S$ is admissible if and only if each summand is of $\mathcal{A}_X$ is of the form $X\wedge F$ or $X\rfloor F$ with $\Delta_{\text{deg}(F)} \Delta_1=-1$. In particular admissibility only depends on the form $F$.
\end{prop}
\begin{prop}\label{propositionadmissible}
If $\mathcal{K}_1\subseteq\big\{\eta\in\Gamma S;\, D^C\eta=0\big\}$. Then for all $\eta,\xi\in\mathcal{K}_1$ the vector field $\{\eta,\xi\}$ is  Killing  if and only if $(D,\mathcal{K}_1)$ is admissible.
\end{prop}

Admissibility of a connection is an algebraic property and has the following consequence.
\begin{lem}\label{Bianchi2}
Let $D$ be an admissible connection  with torsion $\mathcal{T}$. Then  (\ref{BianchiDT}) in Proposition \ref{Bianchi} allows to express the curvature of the Levi-Civita connection in terms of $R$ and $\mathcal{T}$:
\begin{equation}
 R^0_{\kappa\mu\nu\lambda}\gamma^\lambda = ad^C_{R_{\kappa\mu}}\gamma_\nu -\hat D_{[\kappa}\mathcal{T}_{\mu]\nu}\,. 
\end{equation}
\end{lem}


\section{SUSY structures on $M_\q$}

Consider a (pseudo) Riemannian spin manifold $M_\red$ such that $\Delta_1=1$. We use the notations introduced above, i.e.\ $\mathfrak{A}=\Gamma \Lambda S$ with its canonical filtration and canonical deformation $\defo_\q\mathfrak{A}$, the graded manifold $M=(M,\mathfrak{A})$, its deformation $M_\q=(M_\red,\defo_\q\mathfrak{A})$, and the filtration $\mathcal{Z}=\{Z_n\}_{n\in\n}$ such that $\defo_\q^{\mathcal{Z}}\X(M)=\X(M_\q)\cap\X(M)\llbracket\q\rrbracket$.

\begin{defn}\label{defSSKS}
Let $M_\red$ be a (pseudo) Riemannian spin manifold with charge conjugation obeying $\Delta_1=1$. Furthermore let $S$ be a spinor bundle over $M_\red$ and $\mathfrak{A}$ the sheaf of sections in $\Lambda S$. Let $M_\q$ be the deformation of the graded manifold $M=(M_\red,\mathfrak{A})$.  
A {\em SUSY structure}  on  $M_\q$  is a subsuperalgebra  $\mathfrak{S}\subset \defo_\q^\mathcal{Z}\X(M)\subset\X(M_\q)$ which is given by the following data: A connection $D$ on $S$ which determines the splitting (\ref{splitting}),
two $\c$-linear maps
\begin{equation} 
\e:\X(M_\red) \to \X(M_\q)^0\,, \qquad \o:\Gamma S\to \X(M_\q)^1\,, 
\end{equation}
two subspaces $\mathcal{K}_0\subset \X(M_\red) $ and $ \mathcal{K}_1\subset \Gamma S$, as well as a subspace $\mathfrak{Z}\subset \defo_\q^\mathcal{Z}\X(M)$,
such that
\begin{itemize}
\item[a)] $\mathfrak{S}=\e(\mathcal{K}_0)\oplus\o(\mathcal{K}_1)\oplus\mathfrak{Z}$ as superspaces,
\item[b)] $\mathfrak{Z}$  is a module over $\mathfrak{S}\,\text{mod}\,\q^2$,
\item[c)] $\mathfrak{S}\,\text{mod}\,\q$ is a  semidirect product with  non-vanishing commutators
\begin{equation*}
\left.
\begin{array}{rcl}
\big[\e(X),\e(Y)\big] &=&\e([X,Y])\\
\big[\e(X),\o(\eta)\big] &=&\o(\mathcal{L}_X\eta)
\end{array}
\right\}\mod\q
\end{equation*}
for all $X,Y\in\mathcal{K}_0$ and $\eta\in\mathcal{K}_1$, and 
\item[d)] $\mathfrak{S}\,\text{mod}\,\q^2$ is a Lie superalgebra with  non-vanishing commutators 
\begin{equation*}
\left.
\begin{array}{rcl}
	\big[\e(X),\e(Y)\big]	&=&  \e([X,Y])  \\
	\big[\e(X),\o(\eta)\big]	&=& \o( \mathcal{L}_X\eta) \\
	\big[\o(\eta),\o(\xi)\big]	&=& \q \jmath({\{\eta,\xi\}})\\
	\big[\e(X),\q\jmath(Y)\big]	&=& \q\jmath([X,Y])
\end{array}
\right\}\mod\q^2
\end{equation*}
for all $X,Y\in\mathcal{K}_0$ and $\eta,\xi\in\mathcal{K}_1$.
\end{itemize}
The SUSY structure is called {\em mass preserving} if $\e$ and $\o$ preserve the mass dimension.
The elements in $\mathcal{K}_0$ and $\mathcal{K}_1$ are called {\em even} and {\em odd Killing fields}, respectively.
The SUSY structure is called {\em pure} if $\mathfrak{Z}=0\;\text{mod}\,\q^2$ and $\mathfrak{Z}$  is generated by $\mathcal{K}_0$ and $\mathcal{K}_1$, and it is called  {\em finite} if $\mathfrak{Z}$ is finite. 
\end{defn}
If the maps $\o$ and $\e$ preserve the mass dimension as well as the $\z_2$-degree, we see from table \ref{tablemass}  that this is possible only if 
\begin{align*}
\text{image}(\e)&\subset \jmath(\X(M_\red))
		\oplus\Gamma S\otimes\jmath(\Gamma S) 
		\subset Z_0^0 \,, \\
\text{image}(\o)&\subset \jmath(\Gamma S)
		\oplus \q\;\Gamma S\otimes \jmath(\X(M_\red))
		\oplus \q\;\Gamma\Lambda^2S\otimes\jmath(\Gamma S)
		\subset  Z_0^1\oplus  q Z_1^1\,.
\end{align*}

\begin{rem}
\begin{itemize}
\item 
In particular, c) yields that $\mathcal{K}_0$ has to be subalgebra of the vector fields and $\mathcal{K}_1$ has to be invariant under the Lie action of $\mathcal{K}_0$ and $\mathfrak{S}\,\text{mod}\,\q$ is isomorphic to $\mathcal{K}_0\ltimes\mathcal{K}_1$. In d) the even and odd generators are given by $\mathcal{K}_0$ and $\mathcal{K}_1$, respectively.
\item
The parameter which describes the higher order terms of the algebra $\mathfrak{S}$ has mass dimension $1$ and may therefore be interpreted as energy. A useful interpretation of this fact is the following:  For low energies the even fields of infinitesimal transformations  and their odd counterparts decouple but for increasing energies supersymmetric effects have to be taken into account. Moreover, further increasing the energy leads to higher order terms contributing to the algebra.
\item
If $\eta\in\mathcal{K}_1$ obeys $\{\eta,\eta\}=0$, for example $\eta$ is a pure spinor, then $\defo_\q\mathfrak{A}$ is a $Q$-algebra in the sense of \cite{Schwarz-def} by taking $Q:=\o(\eta)$.
\end{itemize}
\end{rem}

\begin{defn}\label{defEO}
Consider the manifolds  $M_\red$, $M$ and $M_\q$ as before. Let $D$ be a spinor connection on $S$. We define 
\begin{align}
\e&:\X(M_\red) \to\X(M_\q)^0 \,,&
\e(X)&:=  \mathcal{L}_X\,,  \\ 
\o&:\Gamma S\to\X(M_\q)^1 \,,&
\o(\varphi)&:=\jmath(\varphi) +\q\;\imath(\varphi)\,.
\end{align}
Here 
\[
\imath:\Gamma S\overset{\textstyle{\imath}_0}{\hookrightarrow}  \Gamma S\oplus \Gamma S_\frac{3}{2}
		=\Gamma S\otimes \X(M_\red)  \overset{\textstyle{\jmath}}{\hookrightarrow} \X(M)^{1,v}
\] 
where $\imath_0$ is the spin-invariant inclusion such that  
$\o(\varphi)(\eta)=\langle\varphi,\eta\rangle +\q\gamma^\mu\varphi\wedge D_\mu\eta$. 
If we write $D=\nabla+\mathcal{A}$,  the image of $\e$ in terms of  $\jmath$ is $\e(X)=\jmath(X)-\nabla X -\mathcal{A}_X$.  
\end{defn}

\begin{thm}\label{purestructure}
We consider the deformation $M_\q$ of  the special graded manifold $M=(M_\red,\mathfrak{A})$. Let $D$ be a connection on the spinor bundle $S$ over the reduced manifold. 
Let $\mathcal{K}_0\subset\X(M_\red) $ be the subset of Killing vector fields on $M_\red$  which leave the connection invariant, i.e. $\mathcal{L}_X D=0$ for all $X\in\mathcal{K}_0$. 
Furthermore let $\mathcal{K}_1\subset\Gamma S$ be a subspace of $D^C$-parallel spinors such that $(D,\mathcal{K}_1)$ is admissible.
Then $\mathcal{K}_0$, $\mathcal{K}_1$ and $\e$, $\o$ cf.\ definition \ref{defEO} define a  mass preserving, pure  SUSY structure on $M_\q$. 
\end{thm}

\begin{rem}
The subset of Killing vector fields that leave the connection $D$ invariant is a subalgebra, due to 
$[L_{[X,Y]},D_Z] = [ L_X, [L_Y,D_Z]]- [ L_Y,[L_X,D_Z]]$.
\end{rem}

\begin{rem}
For a supersymmetry realization similar to  $\mathfrak{S}\,\text{mod}\,\q$ with  $D$ to be the Levi-Civita connection we refer  to \cite{AlCorDevSem}. Moreover, the semidirect product $\mathcal{K}_0\ltimes\mathcal{K}_1$ has also been used in \cite{Habermann3} as an ansatz to construct superalgebras including twistor spinors. The counterexample presented by the author perfectly fits  in our setting, because the Killing connection $D=\nabla+i \gamma$ is not admissible for the used inner product on the spinor bundle.
\end{rem}

In appendix \ref{appenA} we will examine in detail the higher order contributions to the SUSY structure introduced in theorem \ref{purestructure}. This also provides a proof of the statement formulated therein. In particular, we recall that the choice of charge conjugation on the spinor bundle $S$ over the reduced manifold $M_\red$ is in such a way that $\Delta_1=1$. Some of the relations we will use in the following  have also been used  in \cite{Klinker5} when we discussed second order commutators and their  relations to the Bianchi identities cf.\ proposition \ref{Bianchi}.

We will see that all we need  to construct a basis of the  center of the  pure SUSY structure module $\q^{k+1}$ are  mixed commutators of $\ell$ powers of $D$ with $k-\ell$ powers of $\imath$ for $0\leq\ell\leq k$. 
The elementary summands of contributions to the center (modulo $\q^{k+1}$) are shown to  be obtained by the $k$-th power of images of $\imath$ as well as contractions  with images of $\jmath$. We will become more precise at suitable place. Here we will describe the SUSY structure c.f.\ theorem \ref{purestructure} up to terms of order $\q^2$. 

\subsection{The SUSY structure up to order 1}

\begin{lem}\label{lemmapure}
In the situation of theorem \ref{purestructure} we have 
\begin{equation}
 \big[\o(\varphi),\o(\psi)\big] = \q D_{\{\varphi,\psi\}}
      +\q^2  \mathfrak{B}(R;\varphi,\psi) 
      +\q^2 \mathfrak{D}(\mathcal{T};\varphi,\psi) 
	  \label{oo}
\end{equation}
for all $\psi,\varphi\in\mathcal{K}_1$. We use the short notations
\begin{align}
\mathfrak{B}(R;\varphi,\psi) &=\gamma^\mu\varphi\wedge\gamma^\nu\psi\owedge R_{\mu\nu}\,, \label{B}\\
\mathfrak{D}(\mathcal{T};\varphi,\psi) &=
       \frac{1}{2}\big( \gamma^\mu\varphi\wedge  \mathcal{T}_{\mu\nu}\psi 
  	+\gamma^\mu\psi\wedge \mathcal{T}_{\mu\nu} \varphi\big) \otimes D^\nu\,. \label{D}
\end{align}
\end{lem}

For  $X,Y \in\X(M_\red) $, $\varphi,\psi\in \Gamma S$ and  $\Phi\in\Gamma{End}(S)$ the following fundamental commutation relations hold.
\begin{equation}\label{basiccomm}
{\setlength{\arraycolsep}{0.5ex}
\begin{array}{rclrcl}
\multicolumn{6}{c}{\big[\jmath(X),\jmath(Y)\big]		 \, = \, R(X,Y)+\jmath([X,Y])}\,, \\
\big[\jmath(\varphi),\jmath(\psi)\big]		&=& 0\,, 			&\quad
\big[\jmath(X),\jmath(\varphi)\big]		&=& \jmath(D^C_X\phi)\,,\\
 \big[\Phi,\jmath(\varphi)\big]		&=& \jmath(-\Phi^C\varphi)\,,	&\quad
\big[\jmath(X),\Phi\big]			&=& D_X\Phi\,.
\end{array}
}
\end{equation}
From this we get 
\begin{align}
 \big[\jmath(\varphi), \imath(\psi)\big]  = &\  \tfrac{1}{2} D_{\{\varphi,\psi\}}+\gamma^\mu\psi\otimes D^C_\mu\varphi\,,  \label{ji}\\
 \begin{split}
 \big[\imath(\varphi),\imath(\psi)\big] = &\ \gamma^\mu \varphi\wedge \gamma^\nu\psi\owedge R_{\mu\nu} \\
 &    +\gamma^\mu\varphi\wedge \gamma^\nu D^C_\mu\psi \otimes D_\nu
     +\gamma^\mu\psi\wedge \gamma^\nu D^C_\mu\varphi \otimes D_\nu\\
 & 
  +\gamma^\mu\varphi\wedge  ad^C_{\mathcal{A}_\mu}\gamma^\nu\,\psi\otimes D_\nu 
  +\gamma^\mu\psi\wedge ad^C_{\mathcal{A}_\mu}\gamma^\nu\,\varphi\otimes D_\nu\,.
\end{split}
\end{align}
Due to definition \ref{ctorsion} and remark \ref{remarktorsion} the latter yields
\begin{equation}
\begin{split}
\big[\o(\eta),\o(\xi)\big]  =  \q D_{\{\eta,\xi\}}&+\q^2  \gamma^\mu \eta\wedge \gamma^\nu\xi\owedge R_{\mu\nu} \\
			&+\frac{\q^2}{2}\big(\gamma^\mu\eta\wedge \mathcal{T}_{\mu\nu}\xi\otimes D^\nu 
					+\gamma^\mu\xi\wedge \mathcal{T}_{\mu\nu}\eta\otimes D^\nu \big)
\end{split}
\end{equation}
for $\eta,\xi\in\mathcal{K}_1$ which proves lemma \ref{lemmapure}. 
To describe the SUSY structure  we also need 
\begin{equation}
\big[\e(X),D_{\{\eta,\xi\}}\big] 
      =  D_{\{\mathcal{L}_X\varphi,\psi\}} +D_{\{\varphi,\mathcal{L}_X\psi\}} 
\end{equation}
and
\begin{equation} 
\big[\jmath(\varphi) ,D_{\{\eta,\xi\}}\big]  =  0\,,  \label{Dj}
\end{equation}
which holds for all $X\in\mathcal{K}_0$ and $\varphi\in\mathcal{K}_1$.
This is enough to prove the statement of theorem \ref{purestructure} for terms up to order $\q^2$. For the high energy corrections, and therefore for the complete proof, we refer to appendix \ref{appenA}. Before we come to this point we will turn to some examples, in which we will primarily discuss the first order contributions.


\subsection{Examples}

We have the following simple but important example.

\begin{exmp}
Consider flat space $\r^n$ with flat connection on its trivial spinor bundle. Then the SUSY structure c.f.\ theorem \ref{purestructure} admits no higher order corrections and we  recover the usual supersymmetry algebra for $\q=1$.
\end{exmp}
For $D$ to be the Levi-Civita connection we get the next class of examples generalizing the first example.
\begin{exmp}
Consider a Riemannian spin manifold  of dimension $2n$ and holonomy $SU(n)$. The  SUSY structure constructed with $D$ to be the Levi-Civita connection  is finite and $\mathcal{K}_1$  is of dimension two. In particular no torsion terms are present. 
The same but  with $\dim \mathcal{K}_1=1$ holds for Riemannian manifolds of dimension eight or seven and holonomy $Spin(7)$ or $G_2$. The results on supersymmetric Killing structures in \cite{Klinker4} on such manifolds can be obtained by setting $\q=1$. 
\end{exmp}

Now we turn to an example with non vanishing torsion.
\begin{exmp}\label{exampleTorsion}
Given a seven dimensional Riemannian manifold $(N,\hat g)$ and a three dimensional Lorentz manifold $(H,h)$. Consider the product manifold $M=N\times H$  with metric  $g=\hat g \oplus h$ of Lorentzian signature.  The  spinor bundle $S$  of $M$ can be written as the tensor product of the two respective spinor bundles, more precisely $S^+=S_N\otimes S_H$ and the same for $S^-$. The Clifford map decomposes in the same way $\gamma_M=\begin{pmatrix}& \gamma_N\otimes \mathbbm{1}+\mathbbm{1}\otimes \gamma_H \\  \gamma_N\otimes \mathbbm{1}+\mathbbm{1}\otimes \gamma_H & \end{pmatrix}$. We refer to appendix \ref{gammamatrices} for an explicit choice of matrices. In particular this choice is used for explicit calculations.

We specify the above situation and begin with the seven dimensional summand. 
$N$ is assumed to be a manifold which admits a $G_2$ structure $\omega$ parallel with respect to a given metric connection. 
The torsion of this connection shall be contained in the one dimensional $\mathfrak{g}_2$-invariant subspace of $\Lambda^3TN$. Therefore the connection can be written as 
\begin{equation}\label{part1}
D^N_\mu=\nabla_\mu +\tfrac{\lambda}{4}\omega_{\mu\nu\kappa}\gamma^{\nu\kappa}
\end{equation}
and $N$ is an Einstein space with scalar curvature described by the constant $\lambda$ (see \cite{FriedIvanov3}). In this case there exists exactly one spinor $\eta$ such that $D\eta=0$. Furthermore,  the spinor and the invariant form $\omega$ are connected via $\omega \sim C^N_3(\eta,\eta)$. In particular, the only non-trivial forms which can be obtained by $\eta$ are $C^N_3(\eta,\eta)$ and $C^N_0(\eta,\eta)=: |\eta|^2$ as well as their Hodge duals.\footnote{This is due to the fact that $\Delta_0=\Delta_3=\Delta_4=\Delta_7=1$ and $\Delta_1=\Delta_2=\Delta_5=\Delta_6=-1$ on a seven dimensional Riemannian manifold  (compare \cite{Klinker4})} Because $|\eta|^2$ is constant we will consider it to be normed, i.e.\ $|\eta|^2=1$.

The three dimensional part is specified in the following way. 
We consider $H$ to  be a Lie group with totally skew symmetric structure constants $f_{abc}$. 
In particular the Lia algebra in this case is given by $\mathfrak{u}(1)\oplus \mathfrak{u}(1)\oplus \mathfrak{u}(1)$ in the case of vanishing structure constants and $\mathfrak{sl}(2)$ in the remaining case. 
The metric $h$ is given in such a way that it is diagonal with respect to a choice of generators $(E_1,E_2,E_3)$ of the Lie algebra $\mathfrak{h}$ of $H$ with $h(E_1,E_1)=-1,h(E_2,E_2)=h(E_3,E_3)=1$. 
The Levi-Civita connection in this case is determined by $h(\nabla_{E_a}E_b,E_c)=\tfrac{1}{2}f_{abc}$. 
On $H$ and therefore on $S_H$ we consider the flat connection, which in the chosen ON-frame is given by
\begin{equation}\label{part2}
D^H_a =\nabla_a - \tfrac{1}{4}f_{abc}\gamma^{bc}\,.
\end{equation}
Of course all constant spinors are parallel with respect  to this connection. 

On the ten dimensional Lorentzian manifold $M=N\times H$  and also on its  spinor bundle $S$  we consider the connection $D$ given by (\ref{part1}) and (\ref{part2}). 
Furthermore we define $\mathcal{K}_1$ as the (two dimensional) space of $D$-parallel spinors of positive chirality,
\begin{equation}
\mathcal{K}_1= \text{span}\{ \eta_1,\eta_2\},\quad \text{with }\eta_1:=\eta\otimes (1,0)^T, \eta_2=\eta \otimes (0,1)^T\,.
\end{equation}
This turns this example into a candidate of type $IIB$-gravity background. 
The spinors in $\mathcal{K}_1$ yield Killing vector fields
\begin{equation}
X_{\alpha\beta}:=C_1(\eta_\alpha,\eta_\beta)= \tilde\gamma^a_{\alpha\beta}E_a\,,
\end{equation}
or explicitly $X_{11}=i(E_1-E_3)$, $ X_{22}=i(E_1+E_3)$, and $X_{12}=i E_2$. Due to the fact that $|\eta|^2$ is constant with respect to $D$, the $X_{\alpha\beta}$  belong to the set $\mathcal{K}_0$ of Killing vector fields on $M$ that leave the connection invariant. 
The result is summarized as follows\footnote{By $\{e_\alpha\}$ we denote the standard basis of flat $\r^3$ and and by $\{e_\mu,e_8\}$ the basis of $S_N$ induced by the $\mathfrak{g}_2$-structure. Moreover we write $e_{\mu\alpha}=e_\mu\otimes e_\alpha$.} 
\begin{prop}
We consider the manifold $(M,g)$ and connection $D$ as defined above. 
Let $\mathfrak{Z}_2$ be  the lowest order corrections to the SUSY structure cf.\ theorem \ref{purestructure}. A basis of $\mathfrak{Z}_2$ is given by
\begin{align}
\mathfrak{D}(\mathcal{T};\eta_\alpha,\eta_{\alpha'}) 
	& = \eta_1\wedge \eta_2 \otimes D_{X_{\alpha\alpha'}} 
	      +\lambda (\omega_\mu\otimes (e_\alpha\vee e_{\alpha'}) )\otimes D^\mu  
	\label{DDD} 
\end{align}
and $\mathfrak{B}(R;\eta_\alpha,\eta_{\alpha'})$ which acts on a superfunction $\xi\otimes e_\beta\in\Gamma S\subset\Gamma\Lambda S$ by 
\begin{equation}
\mathfrak{B}(R;\eta_\alpha,\eta_{\alpha'})(\xi\otimes e_\beta)
	=2 \xi_\mu R^{\mu\nu\kappa\lambda}e_{\kappa\alpha}\wedge e_{\lambda\alpha'}\wedge e_{\lambda\beta}\,. \label{BBB}
\end{equation} 
\end{prop}

Recalling (\ref{B}) and (\ref{D}), the first order corrections to the SUSY algebra are
\begin{align}
\mathfrak{D}(\mathcal{T};\eta_\alpha,\eta_\beta)& = 
	\mathfrak{D}_\mu(\mathcal{T};\eta_\alpha,\eta_\beta)\otimes D^\mu 
	+ \mathfrak{D}_a(\mathcal{T};\eta_\alpha,\eta_\beta)\otimes D^a\label{D-term} \\
\mathfrak{B}(R;\eta_\alpha,\eta_\beta)&= 
 	\mathfrak{B}_{\mu\nu}(R;\eta_\alpha,\eta_\beta)\owedge\gamma^{\mu\nu} 
	+ \mathfrak{B}_{ab}(R;\eta_\alpha,\eta_\beta)\owedge \gamma^{ab} \label{B-term}
\end{align}
Because $D$ is chosen to be flat in the three dimensional component of $M$, the second summand in (\ref{B-term}) vanishes.
The second summand of (\ref{D-term}) is
\begin{equation}\label{D1}
\mathfrak{D}_c(\mathcal{T};\eta_\alpha,\eta_\beta)
=	f_{abc}\gamma^a\eta_\alpha\wedge \gamma^b\eta_\beta \\
=	\tfrac{1}{2}(\eta\vee\eta) \otimes  \big( f_{abc} \gamma^ae_\alpha \wedge \gamma^be_\beta \big)\,.
\end{equation}
Due to the Fierz identity, which explicitly reads as
\begin{equation}
\eta\otimes\xi=\sum\nolimits_n \frac{\Delta_0(-\Delta_0\Delta_1)^n}{n! \,\text{dim} S}C(\eta,\gamma^{(n)}\xi)\,C\gamma_{(n)}\,,
\end{equation}
and the symmetries $\Delta^3_0=-\Delta^3_1=-1$ we have
\begin{align*}
\gamma^{[a}e_\alpha \wedge\gamma^{b]} e_\beta
=\ &	\frac{1}{2}C^H(\gamma^{[a}e_\alpha, \gamma^{b]}e_\beta) (C^H)^{\alpha\beta}e_\alpha\wedge e_\beta\\
=\ &	 C^H(\gamma^{ab}e_\beta, e_\alpha) e_1\wedge e_2\\
=\ &	 f^{abc} (\tilde\gamma_c)_{\beta\alpha} e_1\wedge e_2\,.
\end{align*}
This together with (\ref{D1}) yields the vertical part of (\ref{DDD}).

If we use the explicit form of the $\gamma$-matrices from appendix \ref{gammamatrices} we see, 
that in these coordinates $\eta=e_8$, i.e.\ $\gamma_\mu\eta=e_\mu$. Therefore 
$\gamma_{[\nu}\eta\wedge  \gamma_{\kappa]}\eta = e_{[\nu}\wedge e_{\kappa]}$
such that the horizontal component of is given by
\begin{align*}
\mathfrak{D}_\mu(\mathcal{T};\eta_\alpha,\eta_\beta)
=\ &	\lambda\omega_{\mu\nu\kappa}\gamma^\nu\eta_\alpha\wedge \gamma^\kappa\eta_\beta \\
=\ &	\tfrac{\lambda}{2}\omega_{\mu\nu\kappa}(e_\nu \wedge e_\kappa)\otimes(e_\alpha\vee e_\beta) \\
=\ & 	\lambda\omega_{\mu}\otimes(e_\alpha\vee e_\beta)\,.
\end{align*}
which is (\ref{DDD}). The second summand in (\ref{B-term}) is 
\begin{align*}
\mathfrak{B}(R;\eta_\alpha,\eta_\beta)
=\ & R^{\mu\nu\kappa\lambda}\gamma_\kappa\eta_\alpha\wedge \gamma_\lambda\eta_\beta \owedge \gamma_{\mu\nu}\\
=\ & \tfrac{1}{2}R^{\mu\nu\kappa\lambda}
	\big( (e_\kappa\wedge e_\lambda)\otimes (e_\alpha\vee e_\beta)\big) \owedge \gamma_{\mu\nu}\\
=\ & \tfrac{1}{2}R^{\mu\nu\kappa\lambda} e_{\kappa\alpha} \wedge e_{\lambda\beta} \owedge \gamma_{\mu\nu}\,.
\end{align*}
We recall that  $D$ is a $\mathfrak{g}_2$-connection and therefore the curvature takes its values in $\mathfrak{g}_2$, i.e. 
$R_{\mu\nu\kappa\lambda}\omega^{\kappa\lambda\sigma}=0$ or equivalently 
$R_{\mu\nu\kappa\lambda}=\tfrac{2}{3}(\delta^{\sigma}_\mu\delta_\nu^\tau +\tfrac{1}{4}*\omega_{\mu\nu}{}^{\sigma\tau})R_{\sigma\tau\kappa\lambda}$ 
which is 
$2 R_{\mu\nu\kappa\lambda} = *\omega_{\mu\nu}{}^{\sigma\tau}R_{\sigma\tau\kappa\lambda}$.
For $\pi,\theta\neq8$ we may decompose $(\gamma^{\mu\nu})_{\pi\theta}$ into the  projections $\Pi_+,\Pi_-$ onto the  summands of the decomposition  $\mathfrak{so}(7)=\mathfrak{g}_2\oplus\mathbf{7}$
\[
(\gamma^{\mu\nu})_{\pi\theta} =-4\cdot \tfrac{2}{3}(\delta_{\pi}^{[\mu}\delta^{\nu]}_\theta +\tfrac{1}{4}*\omega^{\mu\nu}{}_{\pi\theta})
			  +2 \cdot \tfrac{1}{3}(\delta^{\pi}_{[\mu}\delta_{\nu]}^\theta -\tfrac{1}{2}*\omega^{\mu\nu}{}_{\pi\theta})\,.
\]
The result (\ref{BBB}) now follows from
\begin{align*}
R_{\mu\nu\kappa\lambda}\gamma^{\kappa\lambda}\xi 
	& =R_{\mu\nu\kappa\lambda}\big(\gamma^{\kappa\lambda})^{\pi\theta}\xi_\theta e_\pi + (\gamma^{\kappa\lambda})^{\pi 8}\xi_8 e_\pi 
		+(\gamma^{\kappa\lambda})^{8\pi}\xi_\pi e_8\big)  \\
 	& =R_{\mu\nu\kappa\lambda}\big(2 (\Pi_-)^{\kappa\lambda\pi\theta}\xi_\theta e_\pi -4 (\Pi_+)^{\kappa\lambda\pi\theta}\xi_\theta e_\pi
		-\omega^{\kappa\lambda\pi}\xi_8e_\pi +\omega^{\kappa\lambda\pi}\xi_\pi e_8\big)    \\\
	&= 4  \xi^\theta R_{\theta\mu\nu\kappa}e^\kappa \,.
\end{align*}
\end{exmp}


\section{Concluding Remarks and Outlook}
In the preceding sections we described supersymmetry as a semiclassical limit of a geometric structure on a deformed supermanifold. 
After providing the necessary tools  we  defined the SUSY structure  in  definition \ref{defSSKS}.  The example cf.\ theorem \ref{purestructure} was discussed in detail. The constructive proof and in particular  lemma \ref{lemmaAll} provides an explicit  list of all possible high energy corrections to $\mathfrak{S}$  in terms of curvature and torsion of the connection $D$ on the spinor bundle, as well as their covariant derivatives. So this list allows calculations  with explicit examples. Due to lemma \ref{lemmaBasis} the amount of elements in $\mathfrak{Z}$ increase rapidly with raising order which may make the calculations lengthy. 
In many cases some of the terms vanish due to the geometry of the underlying manifold, like in example \ref{exampleTorsion} or the discussion of Killing structures in \cite{Klinker4}. 
Another way to avoid correction terms is to break supersymmetry, i.e.\ we restrict to subsets of the admissible set $\mathcal{K}_1$. This is achieved by imposing vanishing conditions on some of the higher order terms. When we discussed brane metrics in  \cite{Klinker5} the vanishing of $\mathfrak{D}(\mathcal{T},\cdot,\cdot)$ restricted the parallel spinors to those contained in the kernel of the basic vector field of the metric. 
In view of future work one might construct  actions of the infinitesimal automorphism $\mathfrak{S}$ of $M_\q$ on sets of fields which allow  defining generalized invariant Lagrangians. This should be done by taking into account natural bundles over $M_\q$ which  contain deformations of fields with arbitrary spin content. In this respect it might also be interesting in what way quasi morphisms of the algebra $\mathfrak{S}$ can be obtained. They are candidates for duality transformations between different models.


\begin{appendix}

\section{Proof of theorem \ref{purestructure}}\label{appenA}


\subsection{High energy corrections to the SUSY structure} 

The contributions to the SUSY structure of order two or higher in $\q$ have been collected in the set $\mathfrak{Z}$. 
For $\mathfrak{Z}$  to be pure it must  be generated by $\mathcal{K}_1$ and $\mathcal{K}_0$. 
 The following construction will help us to describe the  higher order terms of the SUSY structure cf.\ theorem \ref{purestructure}.

Let $D$ and $\imath$ be abbreviations for the maps from $\Gamma S\otimes \Gamma S$ and $\Gamma S$ to $\X(M)$ given by $\eta\otimes\xi\mapsto D_{\{\eta,\xi\}}$ and $\eta\mapsto \imath(\eta)$. For $k\geq 2$ we consider maps $X^k: \bigoplus_{\ell=k}^{2k}\Gamma S^{\otimes\ell}\to \X(M)$ which images are given by all possible commutators of exactly $k$ images of $D$ and $\imath$. For example in the case $k=5$ such map may be $[D,[\imath,[D,[D,\imath]]]]$.  We define 
\[
U_k:={\rm span}\big\{\text{images of }X^k\big\}\,.
\] 
From the construction it is clear that $U_k\subset \mathcal{Z}_{k}$.
\begin{lem}\label{lemmaBasis}
${U}_2  ={\rm span}\big\{X^2_1(\eta_1,\eta_2),X^2_2(\eta_1,\eta_2,\eta_3),X^2_3(\eta_1,\eta_2,\eta_3,\eta_4) \big| \eta_j\in\Gamma S\big\}$ and for $3\leq k\leq\dim S$ the space $U_k$ is spanned by at most $2^{k-1}$ linear independent types of elements of which $a^k_j$ are $(2k-j)$-linear. Here $a^k_j=\binom{k-3}{j-3}+\binom{k-3}{j-2}+\binom{k-3}{j-1}+\binom{k-3}{j}$ for $j=0,\ldots,\big[\frac{k+2}{2}\big]-1$ and $a^k_j=a^k_{k-j}$ for $\big[\frac{k+2}{2}\big]\leq j\leq k$. If $k>\dim S$ the maximal number of elements is reduced to $\sum_{j=0}^{k-\dim S}a^k_j$.
\end{lem}

\begin{proof}
A basis of $U_k$ will be constructed recursively. Starting with $k=2$  we will for $k\geq3$ get 
$U_k={\rm span}\big\{\text{images of } X^k_1,\ldots, X^k_{2^{k-1}}\big\}$ and, in particular, $a^k_j$ of these maps are $(2k-j)$-linear. 

We start with $U_2$. All maps we get by commutators of $D$ and $\imath$ are 
\begin{align*}
 X_1^2(\eta_1,\eta_2) &= \big[ \imath(\eta_1),\imath(\eta_2) \big] \\
 X_2^2(\eta_1,\eta_2,\eta_3)& = \big[ D_{\{\eta_1,\eta_2\}},\imath(\eta_3)\big] \\
 X_3^2(\eta_1,\eta_2,\eta_3,\eta_4)& = \big[ D_{\{\eta_1,\eta_2\}},D_{\{\eta_3,\eta_4\}} \big] 
\end{align*}
This list  -- with $\eta_j$ varying over $\Gamma S$ -- provides  a basis for $U_2$. 

We list some elements from $U_3$:
\begin{equation*}\label{3basis}
\begin{split}
X^3_1(\eta_1,\ldots,\eta_3)&=\big[ \imath(\eta_1) , \big[\imath(\eta_2)  ,\imath(\eta_3) \big] \big]  \\
X^3_2(\eta_1,\ldots,\eta_4)&=\big[ \imath(\eta_1),\big[ \imath(\eta_2) , D_{\{\eta_3,\eta_4\}} \big] \big] \\
X^3_3(\eta_1,\ldots,\eta_5)&=\big[ D_{\{\eta_1,\eta_2\}}, \big[ D_{\{\eta_3,\eta_4\}},\imath(\eta_5) \big]\big]\\
X^3_4(\eta_1,\ldots,\eta_6)&=\big[ D_{\{\eta_1,\eta_2\}}, \big[ D_{\{\eta_3,\eta_4\}}, D_{\{\eta_5,\eta_6\}} \big]\big]
\end{split}
\end{equation*}
This list -- with $\eta_j$ varying over $\Gamma S$ -- provides a basis for $U_3$  because the missing combinations are given by $\big[D,[\imath,\imath\big]]=\big[\imath,[D,\imath]\big]+\big[\imath,[D,\imath]\big]$ and $\big[\imath,[D,D]\big]=\big[D,[\imath,D]\big]-\big[D,[D,\imath]\big]$ due to the Jacobi identity.

The basis for $U_3$  is the starting point of the recursion process. 
The commutators of the basis of $U_k$ with images of  $\imath$ or $D$ yields elements in $U_{k+1}$. This gives an $(m+1)$-linear or an $(m+2)$-linear element from an $m$-linear element. 

In fact this process yields a basis of $U_k$ because all symmetries which come from the Jacobi identity and which may reduce the dimension have been used in the step from $k=2$ to $k=3$.  In particular more complicated commutators are ruled out by the Jacobi identity as well, e.g.\  for $k=4$ we have $\big[[A,B],[C,D]\big]= \big[A,\big[B,[C,D]\big]\big]- \big[B,\big[A,[C,D]\big]\big]$. 
Therefore, for $k\geq 4$ a basis of $U_k$ is constructed  by 
\[
X^k_\alpha=\big[\imath(\cdot),X^{k-1}_\alpha\big]
\quad\text{and}\quad 
X^k_{2^{k-2}+\alpha}=\big[D_{\{\cdot,\cdot\}},X^{k-1}_\alpha\big]
\] 
for $1\leq\alpha\leq2^{k-2}$.

We denote the maximal number of $(2k-j)$-linear elements in the basis  of $U_k$  by $a_j^k$. The statement on these  numbers  is a consequence of the recursion process and we may deduce them from a combination of four Pascal triangles cf.\ table \ref{numbers}.
\begin{table}[htb]\caption{The maximal number  $a_j^k$}\label{numbers}
$
{\renewcommand{\arraystretch}{1.5}
\begin{array}{cccccccccc}
&k=&3&4&5&6&7&8&9&\ldots \\
&&\cdot&\cdot&\cdot&\cdot&\cdot&\cdot&1&\\
&&\cdot&\cdot&\cdot&\cdot&\cdot&1&&\\\
&&\cdot&\cdot&\cdot&\cdot&1&&7&\\
&&\cdot&\cdot&\cdot&1&&6&&\\
&&\cdot&\cdot&1&&5&&22&\\
&&\cdot&1&&4&&16&&\\
&&1&&3&&11&&42&\\
&\nearrow&&2&&7&&26&&\\
a^k_0&&1&&4&&15&&56&\\
&\nearrow&&2&&8&&30&&\\
a^k_1&&1&&4&&15&&56&\\
&\nearrow&&2&&7&&26&&\\
a^k_2&&1&&3&&11&&42&\\
&\nearrow&&1&&4&&16&&\\
a^k_3&&\nearrow&&1&&5&&22&\\
&a^k_4&&\nearrow&&1&&6&&\\
&&a^k_5&&\nearrow&&1&&7&\\
&&&a^k_6&&\nearrow&&1&&\\
&&&&a^k_7&&\nearrow&&1&\\
&&&&&a^k_8&&\nearrow&&\\
&&&&&&a^k_9&&&
\end{array}}
$
\end{table}

This yields that the maximal number of different types of basis elements in the case $3\leq k\leq\dim S$ is given by $\sum_{j=0}^k a^k_j=2^{k-1}$. 
The number  of elements reduces  for $k>\dim S$ because as a result of the above process the $(2k-j)$-linear vector fields have coefficients in $\Gamma\Lambda^{j}S\oplus\Gamma\Lambda^{j+1}S$ which vanishes for $j=\dim S+1, \ldots, k$. So the corrected value is $\sum_{j=0}^{\dim S}a^k_j$ as stated.
\end{proof}

\begin{prop}\label{collect}
We consider the situation of theorem \ref{purestructure} and expand $\mathfrak{Z}$ as $\mathfrak{Z}=\bigoplus\q^{k}\mathfrak{Z}_k$ with $\mathfrak{Z}_k\subset Z_k$.  Then $\mathfrak{Z}_k=U_k$ where the maps $X^\ell_\alpha$ are restricted to $\mathcal{K}_1$. 
The mass dimension of a $(2k-j)$-linear element of the basis is $\frac{2k-j}{2}$.
\end{prop}

\begin{rem}
The action of $\mathcal{K}_0$ has no effect on the shape of the elements in $U:= \bigoplus_k U_k$ if we restrict the maps to $\mathcal{K}_1$ due to the compatibility of the connection with the infinitesimal transformations of $M_\red$.
\end{rem}

\begin{proof}{[Proposition \ref{collect}]}
The commutators $\big[\o(\mathcal{K}_1),\o(\mathcal{K}_1)\big]$ form a subset of $\q \jmath(\mathcal{K}_0)\oplus\q^2U_2\big|_{\mathcal{K}_1}$. 
The construction of $U$  yields that the restriction of $U_k$ is a subset of $\mathfrak{Z}_k$. 
Due to the Jacobi identity we have $\big[U_k,U_\ell\big]\subset U_{k+\ell}$ such that for the proof of equality of $\mathfrak{Z}_k$ and $U_k$  we need to show that  $\o(\mathcal{K}_1)$ acts on $\bigoplus \q^k U_k$. This is of course true for the second summand of $\o=\jmath+\imath$. For the first summand we recall that due to the Jacobi identity as well as (\ref{Dj}) and (\ref{ji}) any  contraction of one of the basic maps in $U_k$ by $\xi\in\mathcal{K}_1$ can be expressed  as a linear combination of maps from $U_k$. 

A $(2k-j)$-linear basic element is a summand of the result of $(k-j)$ contractions of the $k$-th power of $\imath$, i.e. $[\jmath(\xi_1),[\ldots,[\jmath(\xi_{k-j}),[\imath(\eta_1),[\ldots[\imath(\eta_{k-1}),\imath(\eta_k)]\ldots]]]\ldots]]$, which belongs to $\Lambda^j\Gamma S\otimes\jmath(\X(M_\red))\oplus\Lambda^{j+1}S\otimes\jmath(\Gamma S)$. Taking into account  the prefactor  $\q^{k}$ the mass dimension is $k-\frac{j}{2}=\frac{2k-j}{2}$.
\end{proof}

\begin{lem}\label{lemmaAll}
Each element in $\mathfrak{Z}_k$ is a linear combinations of terms of the form
\begin{equation}\label{AA}
	\mathcal{C}_{\vec\nu_1\ldots\vec\nu_k}\Phi_{(1)}^{\mu_1}{}^{\vec\nu_1}\eta_1
\wedge 	\Phi_{(2)}^{\mu_2}{}^{\vec\nu_2}\eta_2
\wedge	\cdots 
\wedge	\Phi_{(k)}^{\mu_k}{}^{\vec\nu_k}\eta_k
\owedge 	\Psi_{\mu_1\ldots\mu_k}\,,
\end{equation}
and terms of the form
\begin{equation}\label{BB}
	\mathcal{C}_{\vec\kappa_1\ldots\vec\kappa_k}\Omega_{(1)}^{\; \mu}{}^{\vec\kappa_1}\eta_1
\wedge 	\Omega_{(2)}{}^{\vec\kappa_2}\eta_2
\wedge	\cdots 
\wedge	\Omega_{(k)}{}^{\vec\kappa_k} \eta_k
\otimes 	D_\mu\,
\end{equation}
with $\eta_j\in\mathcal{K}_1$ for all $j$, as well as terms with one or more of the $k$ factors contracted by elements in $\mathcal{K}_1$. The contractions $\mathcal{C}$ are in such a way that two factors have at most one index in common.
(\ref{AA}) and (\ref{BB}) obey
\begin{itemize}
\item 
 The length $|\vec\mu|$ of the multi-index $\vec\mu=(\mu_1,\ldots,\mu_k)$ is at least two and at most $k$, and all lengths appear in some summand.
\\
 $\Psi$ is a covariant derivative of the curvature of at most degree $(k-2)$. 
\item
The length of the  multi-index $(\mu_j,\vec\nu_j)$ is at least one  for all $j$ and the length of $(\vec\nu_1,\ldots,\vec\nu_k)$ is $2(k-|\vec\mu|)$. 
\\
$\Phi_{(j)}$ is either  of the form $\hat D^\ell\mathcal{T}$ with  $-1\leq\ell\leq k-3 $ or of the form $ad^C_{D^qR}\hat D^p\mathcal{T}$ with $-1\leq q+p\leq k-4$, $p\geq -1, q\geq 0$. 
\item
The length of the multi-index $(\vec\kappa_j)$ is at least one and the length of $(\vec\kappa_1,\ldots,\vec\kappa_k)$ is $2(k-1)$.
\\
$\Omega_{(j)}$ is either of the form  $\hat D^\ell\mathcal{T}$ of degree $-1\leq\ell\leq k-2 $ or of the form $ad^C_{D^q R}\hat D^p\mathcal{T}$ of at most degree $-1\leq p+q\leq k-3$, $p\geq-1, q\geq0$.
\item   All powers of derivatives from $-1$ up to the stated limits appear in some summand and in each summand the total amount of derivatives is $0$ -- here we naturally define $\hat D_k^{-1}\mathcal{T}:= \gamma_k$  and  the curvature gives a value of two.  
\end{itemize}

\end{lem}

\begin{proof}
We consider the basis of  $\mathfrak{Z}_k$. To get the special type of the summands of each element in $\mathfrak{Z}_k$ we recall that each summand which appears in a commutator of the form $\big[D_{\{\cdot,\cdot\}},X^{k-1}_{j}\big]$ also appears in contractions $\big[\jmath(\cdot),\big[\imath(\cdot),X^{k-1}_{j}\big]\big]$. Therefore, each summand appears in the $k$-th power of $\imath$, i.e.\ in $\big[\imath(\eta_1),\big[\cdots,\big[\imath(\eta_{k-1}),\imath(\eta_k)\big]\cdots\big]\big]$, or one of its contractions. The statement on the shape of the summands is true for $k=2$ and $k=3$ (see below). The step of induction is done by showing that  $\big[\imath(\cdot),Y(\cdot_{k-1})\big]=\sum Y(\cdot_k)$ with $Y(\cdot)$ as in (\ref{AA}) or (\ref{BB}). The fact that all summands appear follows with (\ref{compDad}). 
\end{proof}

\subsection{The lower order basic elements}

We give the explicit form of the images of the maps $X_\alpha^k$ for $k=2,3$ when restricted to $D^C$-parallel, admissible spinors. This yields a description of the basic elements of $\mathfrak{Z}_k$

\underline{$k=2$}\ 
\begin{align*}
 X_1^2(\eta_1,\eta_2) 
	& =	\gamma^\mu\eta_1\wedge\gamma^\nu\eta_2\owedge R_{\mu\nu} \\
&	 \quad	+\tfrac{1}{2}\gamma^\mu\eta_1\wedge \mathcal{T}_{\mu\nu}\eta_2\otimes D^\nu 
		+\tfrac{1}{2}\gamma^\mu\eta_2\wedge \mathcal{T}_{\mu\nu}\eta_1\otimes D^\nu \displaybreak[0] \\
 X_2^2(\eta_1,\eta_2,\eta_3) 
	& = \{\eta_1,\eta_2\}^\mu \gamma^\nu\eta_3\owedge R_{\mu\nu} \\
	&\quad	+\tfrac{1}{2}\{\eta_1,\eta_2\}^\mu \mathcal{T}_{\mu\nu}\eta_3\otimes D^\nu  
 		-\langle \eta_1,\mathcal{T}_{\mu\nu}\eta_2 \rangle \gamma^\mu\eta_3\otimes D^\nu \displaybreak[0]\\
 X_3^2(\eta_1,\ldots,\eta_4)
	& = R(\{\eta_1,\eta_2\},\{\eta_3,\eta_4\}) \\
	&\quad	+\{\eta_1,\eta_2\}^\mu\langle \eta_3,\mathcal{T}_{\mu\nu}\eta_4\rangle D^\nu
		-\{\eta_3,\eta_4\}^\mu\langle \eta_1,\mathcal{T}_{\mu\nu}\eta_2  \rangle D^\nu
\end{align*}
\underline{$k=3$}\ 
\begin{align*}
X^3_1(\eta_1&,\eta_2,\eta_3)=
		 \gamma^\kappa\eta_1\wedge\gamma^\mu\eta_2\wedge\gamma^\nu\eta_3\owedge \big(D_\kappa R\big)_{\mu\nu} \\ 
  	&	
		+\tfrac{1}{2}\Big(  \gamma^\mu\eta_1\wedge\gamma_\kappa\eta_2\wedge\mathcal{T}^{\kappa\nu}\eta_3 
            		+ \gamma^\mu\eta_1\wedge\gamma_\kappa\eta_3\wedge\mathcal{T}^{\kappa\nu}\eta_2 \\
   	&\quad	
		+\gamma_\kappa\eta_1\wedge\gamma^\mu\eta_3 \wedge \mathcal{T}^{\kappa\nu}\eta_2
          		+ \gamma_\kappa\eta_1\wedge \gamma^\mu\eta_2 \wedge \mathcal{T}^{\kappa\nu}\eta_3 \Big) \owedge R_{\mu\nu}\displaybreak[0] \\
  	&	
		+ \tfrac{1}{4}\Big( \gamma^\kappa\eta_1\wedge\mathcal{T}_{\kappa\mu}\eta_2\wedge\mathcal{T}^{\mu\nu}\eta_3
          		+\gamma^\kappa\eta_1\wedge\mathcal{T}_{\kappa\mu}\eta_3\wedge\mathcal{T}^{\mu\nu}\eta_2 \\
   	&\quad	
		+ \mathcal{T}^{\nu\kappa}\eta_1\wedge\gamma^\mu\eta_2\wedge\mathcal{T}_{\mu\kappa}\eta_3 
          		+\mathcal{T}^{\nu\kappa}\eta_1\wedge\gamma^\mu\eta_3\wedge\mathcal{T}_{\mu\kappa}\eta_2  \Big)\otimes D_\nu \\
   	& 
		+\Big(  \tfrac{1}{2}\gamma^\kappa\eta_1\wedge\gamma^\mu\eta_2\wedge \hat D_\kappa\mathcal{T}_{\mu\nu}\eta_3 
         		+ \tfrac{1}{2}\gamma^\kappa\eta_1\wedge\gamma^\mu\eta_3\wedge \hat D_\kappa\mathcal{T}_{\mu\nu}\eta_2 \\
	&\quad 
		- \gamma^\kappa\eta_2\wedge\gamma^\mu\eta_3\wedge ad^C_{R_{\kappa\mu}}\gamma_\nu\eta_1\Big) \otimes D^\nu 
\displaybreak[0]\\
X^3_2(\eta_1&,\ldots,\eta_4)= 
	\{\eta_3,\eta_4\}^\mu\gamma^\kappa\eta_1\wedge \gamma^\nu\eta_2\owedge D_{\kappa}R_{\mu\nu} \\
	&	+\Big( \langle \eta_3,\mathcal{T}^{\kappa\mu}\eta_4\rangle \gamma_\kappa\eta_1\wedge \gamma^\nu\eta_2
		+\langle\eta_3,\mathcal{T}^{\kappa\mu}\eta_4\rangle\gamma^\nu\eta_1\wedge\gamma_\kappa\eta_2\\
	&\quad	
		+\tfrac{1}{2}\{\eta_3,\eta_4\}_\kappa \gamma^\mu\eta_1\wedge \mathcal{T}^{\kappa\nu}\eta_2
		+\tfrac{1}{2}\{\eta_3,\eta_4\}^\mu\gamma_\kappa\eta_1\wedge\mathcal{T}^{\kappa\nu}\eta_2\Big)\owedge R_{\mu\nu}\displaybreak[0] \\
	&	+\tfrac{1}{2}\Big(\langle\eta_3,\mathcal{T}_{\mu\kappa}\eta_4\rangle\gamma^\mu\eta_1\wedge\mathcal{T}^{\kappa\nu}\eta_2
		-\langle\eta_3,\mathcal{T}^{\kappa\nu}\eta_4\rangle\gamma^\mu\eta_1\wedge\mathcal{T}_{\mu\kappa}\eta_2 \\
	&\quad
		+\langle\eta_3,\mathcal{T}_{\mu\kappa}\eta_4\rangle \mathcal{T}^{\nu\kappa}\eta_1\wedge\gamma^\mu\eta_2
		+\tfrac{1}{2}\{\eta_3,\eta_4\}^\mu\mathcal{T}_{\kappa\mu}\eta_2\wedge\mathcal{T}^{\nu\kappa}\eta_1\Big)\otimes D_\nu\\
	& 	+\Big(  \tfrac{1}{2}\{\eta_3,\eta_4\}^\kappa\gamma^\mu\eta_1\wedge\hat D_\mu\mathcal{T}_{\kappa\nu}\eta_2
 		- \langle \eta_3,\hat D_\mu\mathcal{T}_{\kappa\nu}\eta_4\rangle\gamma^\mu\eta_1\wedge \gamma^\kappa\eta_2 \\
	&\quad
		+\{\eta_3,\eta_4\}^\mu\gamma^\kappa\eta_2\wedge ad^C_{R_{\mu\kappa}}\gamma_\nu\eta_1
		\Big)\otimes D^\nu
\displaybreak[0]\\
X^3_3(\eta_1&,\ldots,\eta_5)=
		\{\eta_1,\eta_2\}^\kappa\{\eta_3,\eta_4\}^\mu\gamma^\nu\eta_5\owedge D_\kappa R_{\mu\nu} \\
	&	+\Big( \{\eta_1,\eta_2\}_\kappa\langle\eta_3 ,\mathcal{T}^{\kappa\mu}\eta_4\rangle \gamma^\nu\eta_5 
		- \{\eta_1,\eta_2\}^\mu \langle\eta_3 ,\mathcal{T}^{\kappa\nu}\eta_4\rangle \gamma_\kappa\eta_5 \\
	&\quad 	+ \tfrac{1}{2}\{\eta_1,\eta_2\}_\kappa\{\eta_3,\eta_4\}^\mu\mathcal{T}^{\kappa\nu}\eta_5
		+ \tfrac{1}{2}\{\eta_1,\eta_2\}^\mu\{\eta_3,\eta_4\}_\kappa\mathcal{T}^{\kappa\nu}\eta_5 \Big)\owedge R_{\mu\nu}\displaybreak[0] \\
	&	+\tfrac{1}{2}\Big( \{\eta_1,\eta_2\}^\kappa\langle\eta_3,\mathcal{T}_{\kappa\mu}\eta_4\rangle\mathcal{T}^{\mu\nu}\eta_5
		-\{\eta_1,\eta_2\}^\kappa\langle\eta_3,\mathcal{T}^{\mu\nu}\eta_4\rangle\mathcal{T}_{\kappa\mu}\eta_5 \\
	&\quad 	-\langle\eta_1,\mathcal{T}^{\nu\kappa}\eta_2\rangle\{\eta_3,\eta_4\}^\mu \mathcal{T}_{\mu\kappa}\eta_5
		+\langle\eta_1,\mathcal{T}^{\nu\kappa}\eta_2\rangle \langle\eta_3,\mathcal{T}_{\mu\kappa}\eta_4\rangle\gamma^\mu\eta_5
		\Big)\otimes D_\nu\\
	&	+\Big(\tfrac{1}{2}\{\eta_1,\eta_2\}^\kappa\{\eta_3,\eta_4\}^\mu \hat D_\kappa R_{\mu\nu} 
		-\{\eta_1,\eta_2\}^\kappa\langle\eta_3,\hat D_{\kappa}\mathcal{T}_{\mu\nu}\eta_4\rangle\gamma^\mu\eta_5 \Big)\otimes D^\nu
\displaybreak[0]\\
X^3_4(\eta_1&,\ldots,\eta_6)= 
		\{\eta_1,\eta_2\}^\kappa\{\eta_3,\eta_4\}^\mu\{\eta_5,\eta_6\}^\nu D_\kappa R_{\mu\nu} \\
	& 	+\Big( \{\eta_1,\eta_2\}_\kappa\langle\eta_3,\mathcal{T}^{\kappa\mu}\eta_4\rangle \{\eta_5,\eta_6\}^\nu
		+ \{\eta_1,\eta_2\}^\nu \langle\eta_3,\mathcal{T}^{\kappa\mu}\eta_4\rangle\{\eta_5,\eta_6\}_\kappa \\
	&\quad	- \{\eta_1,\eta_2\}_\kappa \{\eta_3,\eta_4\}^\nu \langle\eta_5,\mathcal{T}^{\kappa\mu}\eta_6\rangle  
		- \{\eta_1,\eta_2\}^\nu \{\eta_3,\eta_4\}_\kappa\langle\eta_5,\mathcal{T}^{\kappa\mu}\eta_6\rangle 
		 \Big) R_{\mu\nu}\displaybreak[0]\\
	&	+\tfrac{1}{2}\Big(
		  \{\eta_1,\eta_2\}^\kappa\langle\eta_3,\mathcal{T}_{\kappa\mu}\eta_4\rangle\langle\eta_5,\mathcal{T}^{\mu\nu}\eta_6\rangle
		-\{\eta_3,\eta_4\}^\kappa\langle\eta_5,\mathcal{T}_{\kappa\mu}\eta_6\rangle\langle\eta_1,\mathcal{T}^{\mu\nu}\eta_2\rangle\\
	&\quad	-\{\eta_1,\eta_2\}^\kappa\langle\eta_5,\mathcal{T}_{\kappa\mu}\eta_6\rangle\langle\eta_3,\mathcal{T}^{\mu\nu}\eta_4\rangle
 		+\{\eta_5,\eta_6\}^\kappa\langle\eta_3,\mathcal{T}_{\kappa\mu}\eta_4\rangle\langle\eta_1,\mathcal{T}^{\mu\nu}\eta_2\rangle
		\Big)D_\nu\\
	& 	+\Big( \{\eta_1,\eta_2\}^\kappa\{\eta_3,\eta_4\}^\mu\langle\eta_5,\hat D_\kappa\mathcal{T}_{\mu\nu}\eta_6\rangle 
		- \{\eta_1,\eta_2\}^\kappa\{\eta_5,\eta_6\}^\mu\langle\eta_3,\hat D_\kappa\mathcal{T}_{\mu\nu}\eta_4\rangle  \Big) D^\nu
\end{align*}

\section{$\gamma$-matrix conventions for $D=3$ and $D=7$}\label{gammamatrices}

For the calculations in example \ref{exampleTorsion} we work with the following explicit $\gamma$-matrix representations. 
On the three dimensional factor we consider $\gamma^H_a$ with   $\gamma^H_1=i\sigma_2$, $\gamma^H_2=\sigma_1$ and $\gamma^H_3=\sigma_3$  as well as the charge conjugation $C^H=\begin{pmatrix}&-1\\1&\end{pmatrix}$.

On the seven dimensional factor we consider -- with respect to an orthonormal frame of $N$  -- the $\gamma$-matrices with
\begin{equation}\label{gamma7-1}
(\gamma_\mu)_{\nu\kappa}= 
	\begin{cases} 
		\omega_{\mu\nu\kappa} 	& \text{if }\mu,\nu\neq8\,,\\  
		\delta_{\mu\nu} 		&\text{if }\kappa=8\,.
	\end{cases}
\end{equation}
They are skew symmetric, such that  $C^N=\mathbbm{1}$,  and they obey 
\begin{equation}\label{gamma7-2}\begin{split}
(\gamma^{\mu\nu})_{\rho\sigma}&=
\begin{cases}
	\omega^{\mu\nu\kappa}(\gamma_\kappa)_{\rho\sigma} -4\delta_{\rho}^{[\mu}\delta^{\nu]}_\sigma 
	=-*\omega^{\mu\nu}{}_{\rho\sigma}-2\delta_{\rho}^{[\mu}\delta^{\nu]}_\sigma
							& \text{if }\rho,\sigma\neq8\,,\\  
	-\omega^{\mu\nu\kappa}(\gamma_\kappa)_{\rho 8}
	= -\omega^{\mu\nu}{}_{\rho}
							& \text{if }\rho\neq \sigma=8\,,
\end{cases}\\
(\gamma_{\mu\nu\kappa})^{\rho\sigma}&=
\begin{cases} 
-\omega_{\mu\nu\kappa}g^{\rho\sigma}-\tfrac{3}{2}\delta^{\{\rho}{}_{[\mu}\omega_{\nu\kappa]}{}^{\sigma\}}
				& \text{if } \rho,\sigma\neq 8\,,\\
*\omega_{\mu\nu\kappa\rho}	& \text{if } \rho\neq \sigma=8\,,\\
-\omega_{\mu\nu\kappa}		& \text{if } \rho=\sigma=8 \,,
\end{cases}
\end{split}\end{equation}
where $\omega$ is the $G_2$-structure of $N$. Some useful  trace identities which are used to get the $\gamma$-matrix identities above are 
\begin{equation}\label{trace}
\begin{aligned}
	\omega^{\mu\nu\kappa}\omega_{\kappa\rho\sigma}
		&= 2\delta^{[\mu}_\rho\delta^{\nu]}_\sigma - *\omega^{\mu\nu}{}_{\rho\sigma}\,,
& 	\omega^{\mu\nu\kappa}\omega_{\rho\nu\kappa}
		&= 6\delta^\mu_\rho\,, \\
	*\omega_{\nu\kappa\rho\lambda}\omega^{\mu\sigma\lambda}
		&=  6\delta^{[\mu}_{[\nu}\omega^{\sigma]}{}_{\kappa\rho]}\,,
&  	*\omega_{\nu\kappa\rho\sigma}\omega_{\mu}{}^{\nu\kappa}
		&=  -4\omega_{\mu\rho\sigma}\,.
\end{aligned}
\end{equation}
We will often consider products of $\gamma$-matrices or there skew symmetrized products with the charge conjugation matrix. The result is denoted by adding a $\tilde{\ }$ to the matrices. 
In particular we have $\tilde\gamma_\mu^N=\gamma^N_\mu$ as well as 
\begin{equation}\label{gamma3-1}
\tilde\gamma^H_1=\begin{pmatrix}-i&\\&-i\end{pmatrix}\,,\ 
\tilde\gamma^H_2=\begin{pmatrix}&i\\i&\end{pmatrix}\,,\ 
\tilde\gamma^H_3=\begin{pmatrix}-i&\\&i\end{pmatrix}\,.
\end{equation}

\end{appendix}

\def\cprime{$'$} \def\cprime{$'$}

\end{document}